\documentclass{article}

\usepackage[a4paper]{geometry}

\usepackage{amsfonts,latexsym,amstext}
\usepackage{amsmath,amscd,amssymb}

\usepackage[latin1]{inputenc}

\usepackage{aeguill}

\usepackage[T1]{fontenc}

\usepackage[usenames]{color}

\usepackage[a4paper]{geometry}

\usepackage{url}

\usepackage{psfrag}

\usepackage{placeins}

\usepackage[dvips]{graphicx}
\usepackage{latexsym,amssymb,amsmath,color,times,bbm}

\newtheorem{thm}{Theorem}[section]
\newtheorem{cor}[thm]{Corollary}
\newtheorem{lem}[thm]{Lemma}
\newtheorem{prop}[thm]{Proposition}
\newtheorem{rem}[thm]{Remark}

\numberwithin{equation}{section}
\newcommand{\norm}[1]{\left\Vert#1\right\Vert}
\newcommand{\abs}[1]{\left\vert#1\right\vert}
\newcommand{\set}[1]{\left\{#1\right\}}

\newcommand{\eps}{\varepsilon}

\newcommand{\esssup}{\mathop{\mbox{ess sup}}}

\newcommand{\essinf}{\mathop{\mbox{ess inf}}}

\def\R{\mathbb R}

\def\Q{\mathbb Q}

\newcommand{\tr}[1]{{\vphantom{#1}}^{\mathit t}{#1}}

\def\sqw{\hbox{\rlap{\leavevmode\raise.3ex\hbox{$\sqcap$}}$%
\sqcup$}}
\def\sqb{\hbox{\hskip5pt\vrule width4pt height6pt depth1.5pt%
\hskip1pt}}

\def\qed{\ifmmode\hbox{\hfill\sqb}\else{\ifhmode\unskip\fi%
\nobreak\hfil
\penalty50\hskip1em\null\nobreak\hfil\sqb
\parfillskip=0pt\finalhyphendemerits=0\endgraf}\fi}
\def\cqfd{\ifmmode\sqw\else{\ifhmode\unskip\fi\nobreak\hfil
\penalty50\hskip1em\null\nobreak\hfil\sqw
\parfillskip=0pt\finalhyphendemerits=0\endgraf}\fi}

\begin{document}

\renewcommand{\labelitemi}{$\bullet$}
\bibliographystyle{plain}
\pagestyle{headings}
\title{On the uniqueness of solutions to quadratic BSDEs with convex generators and unbounded terminal conditions}

\author{
Freddy Delbaen\\
Department of Mathematics\\
ETH-Zentrum, HG G 54.3, CH-8092 Zürich, Switzerland\\
e-mail: delbaen@math.ethz.ch\\ \\
Ying Hu\\
IRMAR, Universit\'e Rennes 1\\ Campus de Beaulieu, 35042 RENNES
Cedex, France\\ e-mail: ying.hu@univ-rennes1.fr\\ \\
Adrien Richou\\
IRMAR, Universit\'e Rennes 1\\ Campus de Beaulieu, 35042 RENNES
Cedex, France\\ e-mail: adrien.richou@univ-rennes1.fr}

\maketitle

\begin{abstract} 
 In \cite{Briand-Hu-08}, the authors proved the uniqueness among the solutions
which admit every exponential moments. In this paper, we prove that
uniqueness
holds among solutions which admit some given exponential moments. These
exponential moments are natural as they are given by the existence theorem. Thanks to this uniqueness result we can strengthen the nonlinear Feynman-Kac formula proved in \cite{Briand-Hu-08}.
\end{abstract}

\section{Introduction}

In this paper, we consider the following quadratic backward stochastic differential equation (BSDE in short for the remaining of the paper)
\begin{equation}
\label{EDSR}
 Y_t=\xi-\int_t^T g(s,Y_s,Z_s)ds +\int_t^T Z_s dW_s, \quad 0 \leqslant t \leqslant T,
\end{equation}
where the generator $-g$ is a continuous real function that is concave and has a quadratic growth with respect to the variable $z$. Moreover $\xi$ is an unbounded random variable (see e.g. \cite{Kobylanski-00} for the case of quadratic BSDEs with bounded terminal conditions). Let us recall that, in the previous equation, we are looking for a pair of process $(Y,Z)$ which is required to be adapted with respect to the filtration generated by the $\mathbb{R}^d$-valued Brownian motion $W$. In \cite{Briand-Hu-08}, the authors prove the uniqueness among the solutions which satisfy for any $p>0$,
$$\mathbb{E}\left[  e^{p \sup_{0 \leqslant t \leqslant T}\abs{Y_t}} \right]<\infty.$$

The main contribution of this paper is to strengthen their uniqueness result. More precisely, we prove the uniqueness among the solutions satisfying:
$$\exists p > \bar{\gamma}, \exists \eps >0, \quad \mathbb{E}\left[e^{p\sup_{0 \leqslant t \leqslant T}\left(Y_t^-+\int_0^t\bar{\alpha}_sds\right)} + e^{\eps\sup_{0 \leqslant t \leqslant T} Y_t^+} \right] < +\infty,$$
where $\bar{\gamma}>0$ and $(\alpha_t)_{t \in [0,T]}$ is a progressively measurable nonnegative stochastic process such that, $\mathbb{P}$-a.s.,
$$\forall (t,y,z) \in [0,T]\times \mathbb{R}\times\mathbb{R}^{1\times d}, \quad g(t,y,z) \leqslant \bar{\alpha}_t+\bar{\beta} \abs{y}+\frac{\bar{\gamma}}{2}\abs{z}^2.$$
Our method is different of that in \cite{Briand-Hu-08} where the authors apply the so-called $\theta$-difference method, i.e. estimating $Y^1-\theta Y^2$, for $\theta\in (0,1)$, and then letting $\theta\rightarrow 0$. Whereas in this paper, we apply a verification method: first we define a stochastic control problem and then we prove that the first component of any solution of the BSDE is the optimal value of this associated control problem. Thus the uniqueness follows immediately. Moreover, using this representation, we are able to give a probabilistic representation of the following PDE:
\begin{equation*}
\partial_t u(t,x)+\mathcal{L}u(t,x)-g(t,x,u(t,x),-\sigma^*\nabla_xu(t,x))=0, \quad u(T,.)=h,
\end{equation*}
where $h$ and $g$ have a ``not too high'' quadratic growth with respect to the variable $x$. We remark that the probabilistic representation is also given by \cite{Briand-Hu-08} under the condition that $h$ and $g$ are subquadratic, i.e.:
$$\forall (t,x,y,z) \in [0,T]\times \mathbb{R}^d \times \mathbb{R}\times\mathbb{R}^{1\times d}, \quad  \abs{h(x)}+\abs{g(t,x,y,z)} \leqslant f(t,y,z)+C\abs{x}^p$$
with $f \geqslant 0$, $C> 0$ and $p<2$.

The paper is organized as follows. In section~2, we prove an existence result in the spirit of \cite{Briand-Hu-06} and \cite{Briand-Hu-08}: here we work with generators $-g$ such that $g^-$ has a linear growth with respect to variables $y$ and $z$. As in part~5 of \cite{Briand-Hu-06}, this assumption allows us to reduce hypothesis of \cite{Briand-Hu-08}. Section~3 is devoted to the optimal control problem from which we get as a byproduct a uniqueness result for quadratic BSDEs with unbounded terminal conditions. Finally, in the last section we derive the nonlinear Feynman-Kac formula in this framework.

Let us close this introduction by giving the notations that we will use in all the paper. For the remaining of the paper, let us fix a nonnegative real number $T>0$. First of all, $(W_t)_{t \in [0,T]}$ is a standard Brownian motion with values in $\mathbb{R}^d$ defined on some complete probability space $(\Omega,\mathcal{F},\mathbb{P})$. $(\mathcal{F}_t)_{t \geqslant 0}$ is the natural filtration of the Brownian motion $W$ augmented by the $\mathbb{P}$-null sets of $\mathcal{F}$. The sigma-field of predictable subsets of $[0,T] \times \Omega$ is denoted $\mathcal{P}$.

As mentioned in the introduction, we will deal only with real valued BSDEs which are equations of type (\ref{EDSR}). The function $-g$ is called the generator and $\xi$ the terminal condition. Let us recall that a generator is a random function $[0,T] \times \Omega \times \mathbb{R} \times \mathbb{R}^{1\times d} \rightarrow \mathbb{R}$ which is measurable with respect to $\mathcal{P} \otimes \mathcal{B}(\mathbb{R}) \otimes \mathcal{B}(\mathbb{R}^{1\times d})$ and a terminal condition is simply a real $\mathcal{F}_T$-measurable random variable. By a solution to the BSDE (\ref{EDSR}) we mean a pair $(Y_t,Z_t)_{t \in [0,T]}$ of predictable processes with values in $\mathbb{R} \times \mathbb{R}^{1\times d}$ such that $\mathbb{P}$-a.s., $t \mapsto Y_t$ is continuous, $t \mapsto Z_t$ belongs to $L^2(0,T)$, $t \mapsto g(t,Y_t,Z_t)$ belongs to $L^1(0,T)$ and $\mathbb{P}$-a.s. $(Y,Z)$ verifies  (\ref{EDSR}). We will sometimes use the notation BSDE($\xi$,$g$) to say that we consider the BSDE whose generator is $g$ and whose terminal condition is $\xi$.

For any real $p\geqslant 1$, $\mathcal{S}^p$ denotes the set of real-valued, adapted and càdlàg processes $(Y_t)_{t \in [0,T]}$ such that
$$\norm{Y}_{\mathcal{S}^p}:=\mathbb{E} \left[\sup_{0 \leqslant t\leqslant T} \abs{Y_t}^p \right]^{1/p} < + \infty.$$
$M^p$ denotes the set of (equivalent class of) predictable processes $(Z_t)_{t \in [0,T]}$ with values in $\mathbb{R}^{1 \times d}$ such that
$$\norm{Z}_{M^p}:=\mathbb{E}\left[\left(\int_0^T \abs{Z_s}^2 ds \right)^{p/2}\right]^{1/p} < +\infty.$$
Finally, we will use the notation $Y^*:=\sup_{0\leqslant t\leqslant T} \abs{Y_t}$ and we recall that $Y$ belongs to the class (D) as soon as the family $\set{Y_{\tau}: \tau\leqslant T \textrm{ stopping time}}$ is uniformly integrable.

\section{An existence result}
In this section, we prove a mere modification of the existence result for quadratic BSDEs obtained in \cite{Briand-Hu-08} by using a method applied in section~5 of \cite{Briand-Hu-08}. We consider here the case where $g^-$ has a linear growth with respect to variables $y$ and $z$. Let us assume the following on the generator.
\paragraph{Assumption (A.1).}
There exist three constants $\beta \geqslant 0$, $\gamma >0$ and $r\geqslant 0$ together with two progressively measurable nonnegative stochastic processes $(\bar{\alpha}_t)_{0\leqslant t \leqslant T}$, $(\underline{\alpha}_t)_{0\leqslant t \leqslant T}$ and a deterministic continuous nondecreasing function $\phi : \mathbb{R}^+ \rightarrow \mathbb{R}^+$ with $\phi(0)=0$ such that, $\mathbb{P}$-a.s.,
\begin{enumerate}
 \item for all $t \in [0,T]$, $(y,z) \mapsto g(t,y,z)$ is continuous;
 \item monotonicity  in $y$: for each $(t,z) \in [0,T] \times \mathbb{R}^{1\times d}$,
$$\forall y \in \mathbb{R}, \quad y(g(t,0,z)-g(t,y,z)) \leqslant \beta \abs{y}^2;$$
 \item growth condition: $\forall (t,y,z) \in [0,T]\times \mathbb{R} \times \mathbb{R}^{1\times d},$
$$-\underline{\alpha}_t -r(\abs{y}+\abs{z}) \leqslant g(t,y,z) \leqslant \bar{\alpha}_t + \phi(\abs{y}) + \frac{\gamma}{2} \abs{z}^2.$$
\end{enumerate}

\begin{thm}
\label{existence}
Let (A.1) hold. 
If there exists $p>1$ such that $$\mathbb{E}\left[\exp \left( \gamma e^{\beta T} \xi^-+\gamma \int_0^T \bar{\alpha}_t e^{\beta t}dt \right) + (\xi^+)^p + \left(\int_0^T \underline{\alpha}_t dt\right)^p\right]<+\infty$$
then the BSDE (\ref{EDSR}) has a solution $(Y,Z)$ such that
$$-\frac{1}{\gamma} \log \mathbb{E} \left[ \exp \left( \gamma e^{\beta(T-t)} \xi^- +\gamma \int_t^T \bar{\alpha}_r e^{\beta(r-t)}dr \right) \Big| \mathcal{F}_t \right] \leqslant Y_t \leqslant Ce^{CT} \left(\mathbb{E} \left[ (\xi^+)^p+\left(\int_t^T \underline{\alpha}_r dr \right)^p \Big| \mathcal{F}_t \right]\right)^{1/p},$$
with $C$ a constant that does not depend on $T$.
\end{thm}

\paragraph{Proof.}
We will fit the proof of Proposition~4 in \cite{Briand-Hu-06} to our situation. Without loss of generality, let us assume that $r$ is an integer. 
For each integer $n\geqslant r$, let us consider the function
$$g_n(t,y,z):=\inf \set{g(t,p,q)+n\abs{p-y}+n\abs{q-z}, (p,q) \in \mathbb{Q}^{1+d}}.$$
$g_n$ is well defined and it is globally Lipschitz continuous with constant $n$. Moreover $(g_n)_{n \geqslant r}$ is increasing and converges pointwise to $g$. Dini's theorem implies that the convergence is also uniform on compact sets. We have also, for all $n\geqslant r$,
$$h(t,y,z):=-\underline{\alpha}_t-r(\abs{y}+\abs{z}) \leqslant g_n(t,y,z) \leqslant g(t,y,z).$$
Let $(Y^n,Z^n)$ be the unique solution in $\mathcal{S}^p \times M^p$ to BSDE($\xi$,$-g_n$). It follows from the classical comparison theorem that
$$ Y^{n+1}_t \leqslant Y^n_t \leqslant Y^r_t.$$
Let us prove that for each $n\geqslant r$
$$Y^n_t \geqslant -\frac{1}{\gamma} \log \mathbb{E} \left[ \exp \left( \gamma e^{\beta(T-t)} \xi^- +\gamma \int_t^T \bar{\alpha}_r e^{\beta(r-t)}dr \right) \Big| \mathcal{F}_t \right]:=X_t.$$
Let $(\tilde{Y}^n,\tilde{Z}^n)$ be the unique solution $\mathcal{S}^p \times M^p$ to BSDE($-\xi^-$,$-g_n^+$). It follows from the classical comparison theorem that $\tilde{Y}^n \leqslant Y^n$ and $\tilde{Y}^n \leqslant 0$. Then, according to Proposition~3 in \cite{Briand-Hu-08}, we have $ \tilde{Y}^n \geqslant X$ and so $ Y^n \geqslant X$ for all $n\geqslant r$.
We set $Y = \inf_{n\geqslant r} Y^n$ and, arguing as in the proof of Proposition~3 in \cite{Briand-Hu-08} or Theorem~2 in \cite{Briand-Hu-06} with a localization argument, we construct a process $Z$ such that $(Y,Z)$ is a solution to BSDE($\xi$,$-g$). For the upper bound, let $(\bar{Y},\bar{Z})$ be the unique solution $\mathcal{S}^p \times M^p$ to BSDE($\xi^+$,$-f$). Then the classical comparison theorem gives us that $Y \leqslant \bar{Y}$ and we apply a classical a priori estimate for $L^p$ solutions of BSDEs in \cite{Briand-Delyon-Hu-Pardoux-Stoica-03} to $\bar{Y}$.\cqfd

\begin{cor}
\label{cor existence}
 Let (A.1) hold. We suppose that $\xi^-+\int_0^T \bar{\alpha}_t dt$ has an exponential moment of order $\gamma e^{\beta T}$ and there exists $p>1$ such that $\xi^+ \in L^p$ and $\int_0^T \underline{\alpha}_t dt \in L^p$. 
\begin{itemize}
 \item If $\xi^- +\int_0^T \bar{\alpha}_tdt$ has an exponential moment of order $qe^{\beta T}$ with $q> \gamma$ then the BSDE (\ref{EDSR}) has a solution $(Y,Z)$ such that $\mathbb{E} \left[ e^{qA^*} \right] < +\infty$ with $A_t := Y_t^-+\int_0^t \bar{\alpha}_s ds$.
 \item If $\xi^++\int_0^T \underline{\alpha}_t dt$ has an exponential moment of order $\eps$ then the BSDE (\ref{EDSR}) has a solution $(Y,Z)$ such that $\mathbb{E} \left[ e^{\eps (Y^+)^*} \right] < +\infty$.
\end{itemize}
\end{cor}

\paragraph{Proof.}
Let us apply the existence result : BSDE (\ref{EDSR}) has a solution $(Y,Z)$ and we have
$$A_t =Y_t^-+\int_0^t \bar{\alpha}_sds \leqslant \frac{1}{\gamma} \log \underbrace{\mathbb{E} \left[ \exp \left( \gamma e^{\beta T} \left(\xi^- +\int_0^T \bar{\alpha}_r dr \right) \right) \Big| \mathcal{F}_t \right]}_{:=M_t}.$$
So $e^{qA_t} \leqslant (M_t)^{q/\gamma}$ with $q/\gamma>1$. Since $M^{p/\gamma}$ is a submartingale, we are able to apply the Doob's maximal inequality to obtain
$$\mathbb{E} \left[ e^{qA^*} \right] \leqslant C_q \mathbb{E} \left[ e^{qe^{\beta T} (\xi^-+\int_0^T\bar{\alpha}_s ds)} \right]<+\infty.$$
To prove the second part of the corollary, we define $$N_t:=\mathbb{E} \left[ (\xi^+)^p+\abs{\int_0^T\underline{\alpha}_s ds}^p \Big| \mathcal{F}_t \right].$$
We set $q>1$. There exists $C_{\eps,p,q}\geqslant 0$ such that $x\mapsto e^{ x^{1/p}\eps/q}$ is convex on $[C_{\eps,p,q}, + \infty[$. We have $e^{\eps/q Y_t^+} \leqslant e^{ (C_{\eps,p,q}+N_t)^{1/p}\eps/q}.$ Since $e^{ (C_{\eps,p,q}+N)^{1/p}\eps/q}$ is a submartingale, we are able to apply the Doob's maximal inequality to obtain
$$ \mathbb{E} \left[ e^{\eps(Y^+)^*} \right] \leqslant C \mathbb{E} \left[ e^{\eps(C_{\eps,p,q}+(\xi^+)^p+(\int_0^T\underline{\alpha}_s ds)^p )^{1/p}} \right] \leqslant C \mathbb{E} \left[ e^{\eps(\xi^++\int_0^T\underline{\alpha}_s ds)} \right]<+\infty.$$
\cqfd

\section{A uniqueness result}
\label{uniqueness result section}
To prove our uniqueness result for the BSDE (\ref{EDSR}), we will introduce a stochastic control problem. For this purpose, we use the following assumption on $g$:
\paragraph{Assumption (A.2).}
There exist three constants $K_{g,y}\geqslant 0$, $\bar{\beta}\geqslant 0$ and $\bar{\gamma} >0$ together with a progressively measurable nonnegative stochastic process $(\bar{\alpha}_t)_{t\in [0,T]}$ such that, $\mathbb{P}$-a.s.,
\begin{itemize}
 \item for each $(t,z) \in [0,T] \times \mathbb{R}^{1 \times d},$
$$\abs{g(t,y,z)-g(t,y',z)}\leqslant K_{g,y} \abs{y-y'}, \quad \forall (y,y') \in \mathbb{R}^2,$$
 \item for each $(t,y,z) \in [0,T]\times  \mathbb{R} \times \mathbb{R}^{1 \times d},$
$$g(t,y,z) \leqslant \bar{\alpha}_t+\bar{\beta}\abs{y}+\frac{\bar{\gamma}}{2}\abs{z}^2,$$
 \item $z\mapsto g(t,y,z)$ is a convex function $\forall (t,y) \in [0,T] \times \mathbb{R}$.
\end{itemize}

Since $g(t,y,.)$ is a convex function we can define the Legendre-Fenchel transformation of $g$ :
$$f(t,y,q)=\sup_{z} \left(zq-g(t,y,z)\right), \quad \forall t \in [0,T], q \in \mathbb{R}^{d}, y \in \mathbb{R}.$$
$f$ is a function with values in $\mathbb{R} \cup \set{+\infty}$ that verifies direct properties.
\begin{prop}
$ $
\begin{itemize}
 \item $\forall (t,y,y',q) \in [0,T]\times \mathbb{R}\times \mathbb{R}\times\mathbb{R}^d$ such that $f(t,y,q)<+\infty$,
$$f(t,y',q)<+\infty \quad \textrm{ and } \quad \abs{f(t,y,z)-f(t,y',z)} \leqslant K_{g,y} \abs{y-y'}.$$ 
 \item $f$ is a convex function in $q$,
 \item $f(t,y,q) \geqslant -\bar{\alpha}_t-\bar{\beta}\abs{y}+\frac{1}{2\bar{\gamma}}\abs{q}^2$.
\end{itemize}
\end{prop}
We set $N \in \mathbb{N}^*$ such that 
\begin{equation}
\label{definition N}
\frac{T}{N} < \left( \frac{1}{\bar{\gamma}} -\frac{1}{p} \right) \frac{1}{\bar{\beta}(1/p+1/\eps)}.
\end{equation}
For $i \in \set{0,...,N}$ we define $t_i:=\frac{iT}{N}$ and
\begin{eqnarray*}
\mathcal{A}_{t_i,t_{i+1}}(\eta)&:=&\left\{(q_s)_{s \in [t_i,t_{i+1}]}, \quad \int_{t_i}^{t_{i+1}} \abs{q_s}^2ds < +\infty \,\,\,\, \mathbb{P}-a.s., \right. \\
& &(M^i_t)_{t \in [t_i,t_{i+1}]} \textrm{ is a martingale}, \quad \mathbb{E}^{\mathbb{Q}^i} \left[\abs{\eta}+\int_{t_i}^{t_{i+1}} \abs{f(s,0,q_s)}ds \right]<+\infty, \\
& & \left. \textrm{ with } M^i_t:=\exp \left(\int_{t_i}^t q_s dW_s-\frac{1}{2}\int_{t_i}^t \abs{q_s}^2ds \right) \textrm{ and } \frac{d\mathbb{Q}^i}{d\mathbb{P}}:=M_{t_{i+1}}^i \right\}.
\end{eqnarray*}
Let $q$ be in $\mathcal{A}_{t_i,t_{i+1}}(\eta)$. We define $dW^q_t:=dW_t-q_tdt$. Thanks to the Girsanov theorem, $(W_{t_i+h}^q-W_{t_i}^q)_{h \in [0,1/N]}$ is a Brownian motion under the probability $\mathbb{Q}^i$. So, we are able to apply Proposition 6.4 in \cite{Briand-Delyon-Hu-Pardoux-Stoica-03} to show this existence result:
\begin{prop}
 There exist two processes $(Y^{\eta,q},Z^{\eta,q})$ such that $(Y^{\eta,q}_t)_{t \in [t_i,t_{i+1}]}$ belongs to the class (D) $\int_{t_i}^{t_{i+1}} \abs{Z_s^{\eta,q}}^2 ds <+ \infty$ $\mathbb{P}-a.s.$, $\int_{t_i}^{t_{i+1}} \abs{f(s,Y_s^{\eta,q},q_s)} ds <+ \infty$ $\mathbb{P}-a.s.$ and
$$Y_t^{\eta,q}=\eta+\int_t^{t_{i+1}} f(s,Y_s^{\eta,q},q_s)ds+\int_t^{t_{i+1}}Z_s^{\eta,q} dW_s^q, \quad t_i \leqslant t \leqslant t_{i+1}.$$
\end{prop}
We are now able to define the admissible control set:
\begin{eqnarray*}
\mathcal{A}&:=&\left\{(q_s)_{s \in [0,T]}, \quad q_{|[t_{N-1},T]} \in \mathcal{A}_{t_{N-1},T}(\xi), \forall i \in \set{N-2,\dots,0}, q_{|[t_i,t_{i+1}]} \in \mathcal{A}_{t_i,t_{i+1}}\left(Y_{t_{i+1}}^{q}\right)\right.\\
& &  \left. \textrm{ with } Y_{t_{i+1}}^q:=Y_{t_{i+1}}^{Y_{t_{i+2}}^q,q_{|[t_{i+1},t_{i+2}]}} \textrm{ and } Y_T^q:=\xi \right\}.
\end{eqnarray*}
$\mathcal{A}$ is well defined by a decreasing recursion on $i \in \set{0,\dots,N-1}$. For $q \in \mathcal{A}$ we can define our cost functional $Y^q$ on $[0,T]$ by
$$\forall i \in \set{N-1,\dots,0}, \forall t \in [t_i,t_{i+1}], \quad Y^q_t:=Y_t^{Y_{t_{i+1}}^q,q_{|[t_i,t_{i+1}]}}.$$
$Y^q$ is also well defined by a decreasing recursion on $i \in \set{0,\dots,N-1}$. Finally, the stochastic control problem consists in minimizing $Y^q$ among all the admissible controls $q \in \mathcal{A}$. Our strategy to prove the uniqueness is to prove that given a solution $(Y,Z)$, the first component is the optimal value.
\begin{thm}
\label{thm unicite}
We suppose that there exists a solution $(Y,Z)$ of the BSDE (\ref{EDSR}) verifying 
$$\exists p > \bar{\gamma}, \exists \eps >0, \quad \mathbb{E}\left[\exp\left(pA^*\right) + \exp\left(\eps(Y^+)^*\right) \right] < +\infty,$$
with $A_t:=Y_t^-+\int_0^t\bar{\alpha}_s ds$. Then we have $Y= \essinf_{q \in \mathcal{A}} Y^q$, and there exists $q^* \in \mathcal{A}$ such that $Y=Y^{q^*}$. Moreover, this implies that the solution $(Y,Z)$ is unique among solutions verifying such condition.
\end{thm}
\paragraph{Proof.}
Let us first prove that for any $q$ admissible, we have $Y \leqslant Y^q$. To do this, we will show that $Y_{|[t_i,t_{i+1}]} \leqslant Y^q_{|[t_i,t_{i+1}]}$ by decreasing recurrence on $i \in \set{0,N-1}$. Firstly, we have $Y_T=Y^q_T=\xi$. Then we suppose that $Y_t \leqslant Y_t^q$, $\forall t \in [t_{i+1},T]$. We set $t \in [t_i,t_{i+1}]$,
$$\tau_n^i:=\inf \set{s\geqslant t, \sup\set{\int_{t}^s \abs{Z_u}^2du,\int_{t}^s \abs{Z_u^q}^2du,\int_{t}^s \abs{q_u}^2du}>n} \wedge t_{i+1},$$
$h(s,y,z):=-g(s,y,z)+zq_s$, and
\begin{displaymath}
 h_s:=\left\{
\begin{array}{cl}
 \dfrac{ h(s,Y_s^q,Z_s)-h(s,Y_s,Z_s) }{ Y_s^q-Y_s } & \textrm{if }Y_s^q-Y_s \neq 0\\
0 &\textrm{otherwise.}\\
\end{array}
\right.
\end{displaymath}
We observe that $\abs{h_s} \leqslant K_{g,y}$. Then, by applying Itô formula to the process $(Y_s^q-Y_s)e^{\int_{t}^s h_udu}$ we obtain
$$Y_t^q-Y_t=e^{\int_t^{\tau_n^i} h_sds}\left[Y_{\tau_n^i}^q-Y_{\tau_n^i}\right]+\int_t^{\tau_n^i} e^{\int_t^{s} h_udu}\left[f(s,Y_s^q,q_s)-h(s,Y_s^q,Z_s)\right]ds+\int_t^{\tau_n^i}e^{\int_t^{s} h_udu}\left[ Z_s^q-Z_s\right] dW_s^q.$$
By definition, $f(s,Y_s^q,q_s)-h(s,Y_s^q,Z_s)\geqslant 0$, so 
$$Y_t^q-Y_t \geqslant \mathbb{E}^{\mathbb{Q}^i} \left[\left.e^{\int_t^{\tau_n^i} h_sds}\left[Y_{\tau_n^i}^q-Y_{\tau_n^i}\right] \right| \mathcal{F}_t \right].$$
Since $\left(Y_{\tau_n^i}^q e^{\int_t^{\tau_n^i}h_s ds} \right)_n$ tends to $Y_{t_{i+1}}^q e^{\int_t^{t_{i+1}}h_s ds}$ almost surely and is uniformly integrable, we have
$$\lim_{n\rightarrow + \infty} \mathbb{E}^{\mathbb{Q}^i} \left[ \left. e^{\int_t^{\tau_n^i} h_sds}Y_{\tau_n^i}^q \right| \mathcal{F}_t \right]= \mathbb{E}^{\mathbb{Q}^i} \left[ \left. e^{\int_t^{t_{i+1}} h_sds}Y_{t_{i+1}}^q \right| \mathcal{F}_t \right].$$
Moreover, $\abs{Y_{\tau_n^i} e^{\int_t^{\tau_n^i} h_s ds}} \leqslant (Y^+)^*e^{TK_{g,y}}+(Y^-)^*e^{TK_{g,y}}$, so, by the dominated convergence theorem we obtain 
$$\lim_{n\rightarrow + \infty} \mathbb{E}^{\mathbb{Q}^i} \left[\left. e^{\int_t^{\tau_n^i} h_sds}Y_{\tau_n^i} \right| \mathcal{F}_t \right]= \mathbb{E}^{\mathbb{Q}^i} \left[\left. e^{\int_t^{t_{i+1}} h_sds}Y_{t_{i+1}} \right| \mathcal{F}_t \right].$$
Finally,
$$Y_t^q-Y_t \geqslant \lim_{n \rightarrow + \infty} \mathbb{E}^{\mathbb{Q}^i} \left[e^{\int_t^{\tau_n^i} h_sds}\left[Y_{\tau_n^i}^q-Y_{\tau_n^i}\right] \left| \mathcal{F}_t \right.\right]=\mathbb{E}^{\mathbb{Q}^i} \left[\left. e^{\int_t^{t_{i+1}} h_sds}\left(Y_{t_{i+1}}^q-Y_{t_{i+1}}\right)\right|\mathcal{F}_t\right] \geqslant 0,$$
because $Y_{t_{i+1}}^q\geqslant Y_{t_{i+1}}$ by the recurrence's hypothesis.

Now we set $\tr{q^*_s} \in \partial_z g(s,Y_s,Z_s)$ with $\partial_z g(s,Y_s,Z_s)$ the subdifferential of $z\mapsto g(s,Y_s,z)$ at $Z_s$. We recall that for a convex function $l:\mathbb{R}^{1 \times d} \rightarrow \mathbb{R}$, the subdifferential of $l$ at $x_0$ is the non-empty convex compact set of $u \in \mathbb{R}^{1 \times d}$ such that 
$$l(x)-l(x_0) \geqslant u\tr{(x-x_0)}, \quad \forall x \in \mathbb{R}^{1 \times d}.$$
We have $f(s,Y_s,q_s^*)=zq^*_s-g(s,Y_s,Z_s)$ for all $s \in [0,T]$, so 
\begin{eqnarray*}
 g(s,Y_s,Z_s) &\leqslant& Z_sq_s^*-\frac{1}{2\bar{\gamma}}\abs{q_s^*}^2+\bar{\beta} \abs{Y_s}+\bar{\alpha}_s\\
&\leqslant& \frac{1}{2}\left( 2\bar{\gamma} \abs{Z_s}^2+\frac{\abs{q_s^*}^2}{2\bar{\gamma}} \right) -\frac{1}{2\bar{\gamma}}\abs{q_s^*}^2 +\bar{\beta} \abs{Y_s}+\bar{\alpha}_s\\
\frac{\abs{q_s^*}^2}{4\bar{\gamma}} &\leqslant& -g(s,Y_s,Z_s)+ \bar{\gamma} \abs{Z_s}^2 + \bar{\beta} \abs{Y_s}+\bar{\alpha}_s,\\
\end{eqnarray*}
and finally, $\int_{0}^{T} \abs{q_s^*}^2ds < +\infty$, $\mathbb{P}$-a.s.. Moreover, $\forall t,t' \in [0,T]$,
$$Y_t=Y_{t'}+\int_t^{t'} f(s,Y_s,q_s^*)ds+\int_t^{t'}Z_s(dW_s+q_s^* ds).$$
Thus, we just have to show that $q^*$ is admissible to prove that $q^*$ is optimal, i.e. $Y=Y^{q^*}$. For this, we must prove that $(q_s^*)_{s \in [t_i,t_{i+1}]} \in \mathcal{A}_{t_i,t_{i+1}}(Y_{t_{i+1}})$ for $i \in \set{0,\dots,N-1}$. We define 
$$M^i_t:=\exp\left(\int_{t_i}^t q_s^*dW_s-\frac{1}{2}\int_{t_i}^t \abs{q_s^*}^2 ds \right), \quad \frac{d\mathbb{Q}^{*,i}}{d\mathbb{P}}:=M^i_{t_{i+1}},$$
$$\tau_n^i=\inf \set{t \in [t_i,t_{i+1}], \sup \left(\int_{t_i}^t \abs{q_s^*}^2 ds,\int_{t_i}^t \abs{Z_s}^2 ds\right) >n}\wedge t_{i+1},  \quad \frac{d\mathbb{Q}^{*,i}_n}{d\mathbb{P}}:=M^i_{\tau_{n}^i}.$$
Let us show the following lemma:
\begin{lem}
\label{lemme martingale uniformement integrable}
 $(M^i_{\tau_n^i})_n$ is uniformly integrable.
\end{lem}
\paragraph{Proof.}
Firstly, from
$$xy \leqslant \exp(x)+y(\log(y)-1),\quad \forall (x,y) \in \mathbb{R}\times \mathbb{R}^{+*},$$
we deduce
\begin{equation}
\label{inegalite exp log}
xy=px\frac{y}{p} \leqslant \exp(px)+\frac{y}{p}\left(\log y -\log p -1\right).
\end{equation}
Thus
\begin{eqnarray*}
\mathbb{E}^{\mathbb{Q}_n^{*,i}} \left[ A^* \right] &=& \mathbb{E}  \left[ M^i_{\tau_n^i} A^* \right] \leqslant \mathbb{E} \left[ \exp (pA^*) \right] +\frac{1}{p} \mathbb{E} \left[ M^i_{\tau_n^i} \left(\log M^i_{\tau_n^i} -\log p -1\right) \right]\\
&\leqslant & C_p+\frac{1}{2p} \mathbb{E}^{\mathbb{Q}_n^{*,i}} \left[ \int_{t_i}^{\tau_n^i} \abs{q_s^*}^2 ds\right],
\end{eqnarray*}
and, in the same manner,
\begin{eqnarray*}
\mathbb{E}^{\mathbb{Q}_n^{*,i}} \left[ (Y^+)^* \right] &\leqslant & C_{\eps}+\frac{1}{2\eps} \mathbb{E}^{\mathbb{Q}_n^{*,i}} \left[ \int_{t_i}^{\tau_n^i} \abs{q_s^*}^2 ds\right].\\
\end{eqnarray*}
Since $g(s,Y_s,Z_s)=Z_s q_s^*-f(s,Y_s,q_s^*)$ and $(M^i_{t \wedge \tau_n^i})_{t \in [t_i,t_{i+1}]}$ is a martingale, we can apply the Girsanov theorem and we obtain
$$\mathbb{E}^{\mathbb{Q}_n^{*,i}} \left[ Y_{\tau_n^i}+\int_{t_i}^{\tau_n^i} f(s,Y_s,q_s^*)ds \right] = \mathbb{E}^{\mathbb{Q}_n^{*,i}} \left[Y_{t_i} \right]= \mathbb{E} \left[M_{\tau_n^i}^i Y_{t_i} \right] = \mathbb{E}\left[Y_{t_i} \right].$$
Moreover $f(t,y,q)\geqslant \frac{1}{2\bar{\gamma}}\abs{q}^2-\bar{\beta}\abs{y}-\bar{\alpha}_t$ and $Y_{\tau_n^i} \geqslant -Y_{\tau_n^i}^-$, so
\begin{eqnarray*}
 \mathbb{E} \left[Y_{t_i} \right] &\geqslant& -\mathbb{E}^{\mathbb{Q}_n^{*,i}} \left[ Y_{\tau_n^i}^- \right]- \mathbb{E}^{\mathbb{Q}_n^{*,i}} \left[\int_{t_i}^{\tau_n^i} \bar{\alpha}_s ds \right] +\frac{1}{2\bar{\gamma}} \mathbb{E}^{\mathbb{Q}_n^{*,i}} \left[\int_{t_i}^{\tau_n^i} \abs{q_s^*}^2ds \right]-\bar{\beta} \mathbb{E}^{\mathbb{Q}_n^{*,i}} \left[\int_{t_i}^{\tau_n^i} \abs{Y_s}ds \right]\\
&\geqslant & C-\mathbb{E}^{\mathbb{Q}_n^{*,i}} \left[ A^* \right] + \frac{1}{2\bar{\gamma}} \mathbb{E}^{\mathbb{Q}_n^{*,i}} \left[ \int_{t_i}^{\tau_n^i} \abs{q_s^*}^2ds \right]
 - \frac{T}{N}\left( \bar{\beta} \mathbb{E}^{\mathbb{Q}_n^{*,i}} \left[ (Y^-)^* + (Y^+)^* \right] \right)\\
&\geqslant & C_{p,\eps}+\frac{1}{2}\underbrace{\left(\frac{1}{\bar{\gamma}}-\frac{1}{p}-\frac{T}{N}\left( \frac{\bar{\beta}}{p}+\frac{\bar{\beta}}{\eps} \right) \right)}_{>0} \mathbb{E}^{\mathbb{Q}_n^{*,i}} \left[ \int_{t_i}^{\tau_n^i} \abs{q_s^*}^2ds \right].
\end{eqnarray*}
This inequality explains why we take N verifying (\ref{definition N}).  Finally we get that 
\begin{equation}
\label{esperance q majoree}
2\mathbb{E} \left[M^i_{\tau_n^i}\log  M^i_{\tau_n^i} \right]=\mathbb{E}^{\mathbb{Q}_n^{*,i}} \left[ \int_{t_i}^{\tau_n^i} \abs{q_s^*}^2ds \right] \leqslant C_{p,\eps}.
\end{equation}
Then we conclude the proof of the lemma by using the de La Vallée Poussin lemma.
\cqfd
Thanks to this lemma, we have that $\mathbb{E}[M_{t_{i+1}}^i]=1$ and so $(M_t^i)_{t \in [t_i,t_{i+1}]}$ is a Martingale. Moreover, applying Fatou's lemma and inequality (\ref{esperance q majoree}), we obtain
\begin{equation}
\label{Fatou}
2\mathbb{E}\left[M^i_{t_{i+1}}\log  M^i_{t_{i+1}} \right]=\mathbb{E}^{\mathbb{Q}^{*,i}} \left[\int_{t_i}^{t_{i+1}} \abs{q_s^*}^2ds \right] \leqslant \liminf_n \mathbb{E}^{\mathbb{Q}_n^{*,i}} \left[\int_{t_i}^{\tau_n^i} \abs{q_s^*}^2ds \right] < +\infty.
\end{equation}
So, by using this result and inequality (\ref{inegalite exp log}) we easily show that $\mathbb{E}^{\mathbb{Q}^{*,i}} \left[ (Y^+)^* +(Y^-)^* \right] < + \infty$.
To conclude we have to prove that $\mathbb{E}^{\mathbb{Q}^{*,i}} \left[\int_{t_i}^{t_{i+1}} \abs{f(s,0,q_s^*)}ds \right]<+\infty$:

\begin{eqnarray*}
\mathbb{E}^{\mathbb{Q}^{*,i}} \left[\int_{t_i}^{t_{i+1}} \abs{f(s,0,q_s^*)}ds \right] &\leqslant & \mathbb{E}^{\mathbb{Q}^{*,i}} \left[\int_{t_i}^{t_{i+1}} \abs{f(s,Y_s,q_s^*)}+K_{g,y}\abs{Y_s}ds \right]\\
&\leqslant & \mathbb{E}^{\mathbb{Q}^{*,i}} \left[\int_{t_i}^{t_{i+1}} \abs{f(s,Y_s,q_s^*)}ds+K_{g,y}T\left( (Y^+)^* +(Y^-)^* \right) \right]\\
& \leqslant& C+ \mathbb{E}^{\mathbb{Q}^{*,i}} \left[\int_{t_i}^{t_{i+1}} f^+(s,Y_s,q_s^*)+f^-(s,Y_s,q_s^*)ds \right].
\end{eqnarray*}
Firstly,
$$\mathbb{E}^{\mathbb{Q}^{*,i}} \left[\int_{t_i}^{t_{i+1}} f^-(s,Y_s,q_s^*)ds \right] \leqslant \mathbb{E}^{\mathbb{Q}^{*,i}} \left[\int_{t_i}^{t_{i+1}} \bar{\alpha}_s+\bar{\beta}\abs{Y_s}ds \right]<+\infty.$$
Moreover, thanks to the Girsanov theorem we have
$$\mathbb{E}^{\mathbb{Q}^{*,i}} \left[Y_{t_i} \right] = \mathbb{E}^{\mathbb{Q}^{*,i}} \left[ Y_{\tau_n^i}+\int_{t_i}^{\tau_n^i} f(s,Y_s,q_s^*)ds \right],$$
so
\begin{eqnarray*}
\mathbb{E}^{\mathbb{Q}^{*,i}} \left[ \int_{t_i}^{\tau_n^i} f^+(s,Y_s,q_s^*)ds\right] &\leqslant & \mathbb{E}^{\mathbb{Q}^{*,i}} \left[ Y_{t_i}-Y_{\tau_n^i} \right] + \mathbb{E}^{\mathbb{Q}^{*,i}} \left[ \int_{t_i}^{\tau_n^i} f^-(s,Y_s,q_s^*)ds\right]\\
&\leqslant & C + \mathbb{E}^{\mathbb{Q}^{*,i}} \left[ \int_{t_i}^{t_{i+1}} f^-(s,Y_s,q_s^*)ds\right] \leqslant C
\end{eqnarray*}
Finally, $\mathbb{E}^{\mathbb{Q}^{*,i}} \left[ \int_{t_i}^{t_{i+1}} f^+(s,Y_s,q_s^*)ds\right]<+\infty$ and
$\mathbb{E}^{\mathbb{Q}^{*,i}} \left[\int_{t_i}^{t_{i+1}} \abs{f(s,0,q_s^*)}ds \right] < +\infty$.
Thus, we prove that $q^*$ is optimal, i.e. $Y^{q^*}=Y$.

The uniqueness of $Y$ is a mere consequence of the fact that $Y=Y^{q^*}=\essinf_{q\in \mathcal{A}} Y^q$. The uniqueness of $Z$ follows immediately.
\cqfd
\begin{rem}
\label{remarque espe q2 bornee}
By taking into consideration the inequality (\ref{Fatou}) it is possible to restrict the admissible control set by considering
$$\tilde{\mathcal{A}}_{t_i,t_{i+1}}(\eta):=\mathcal{A}_{t_i,t_{i+1}}(\eta) \cap \set{(q_s)_{s \in [t_i,t_{i+1}]}, \quad \mathbb{E}^{\mathbb{Q}^{i}} \left[\int_{t_i}^{t_{i+1}} \abs{q_s}^2ds \right] < + \infty}$$
instead of $\mathcal{A}_{t_i,t_{i+1}}(\eta)$.
\end{rem}

\begin{rem}
\label{N=1}
If we have $g(t,y,z) \leqslant g(t,0,z)$, then $f(t,y,q) \geqslant f(t,0,q) \geqslant \frac{1}{2\bar{\gamma}}\abs{q}^2-\bar{\alpha}_t$ and we do not have to introduce $N$ in the proof of lemma~\ref{lemme martingale uniformement integrable}. So we have a simpler representation theorem:
$$Y_t=\essinf_{q \in \mathcal{A}_{0,T}(\xi)} Y_t^q, \quad \forall  t \in [0,T].$$
For example, when $g$ is independent of $y$, we obtain
$$Y_t=\essinf_{q \in \mathcal{A}_{0,T}(\xi)} \mathbb{E}^{\mathbb{Q}}\left[\xi+\int_t^T f(s,q_s)ds \Big| \mathcal{F}_t \right], \quad \forall  t \in [0,T].$$
\end{rem}

\section{Application to quadratic PDEs}
In this section we give an application of our results concerning BSDEs to PDEs which are quadratic with respect to the gradient of the solution. Let us consider the following semilinear PDE
\begin{equation}
 \label{EDP}
\partial_t u(t,x)+\mathcal{L}u(t,x)-g(t,x,u(t,x),-\sigma^*\nabla_xu(t,x))=0, \quad u(T,.)=h,
\end{equation}
where $\mathcal{L}$ is the infinitesimal generator of the diffusion $X^{t,x}$ solution to the SDE
\begin{equation}
 \label{EDS}
X_s^{t,x}=x+\int_{t}^s b(u,X_u^{t,x})ds + \int_{t}^s \sigma(u) dW_u, \quad t\leqslant s \leqslant T, \textrm{ and } X_s^{t,x}=x, \quad s\leqslant t.
\end{equation}
The nonlinear Feynman-Kac formula consists in proving that the function defined by the formula 
\begin{equation}
 \label{feynman-kac}
\forall (t,x) \in [0,T] \times \mathbb{R}^d, \quad u(t,x):=Y_t^{t,x}
\end{equation}
where, for each $(t_0,x_0) \in [0,T] \times \mathbb{R}^d$, $(Y^{t_0,x_0},Z^{t_0,x_0})$ stands for the solution to the following BSDE
\begin{equation}
 \label{EDSR markovienne}
Y_t=h(X_T^{t_0,x_0})-\int_t^T g(s,X_s^{t_0,x_0},Y_s,Z_s)ds-\int_t^T Z_s dW_s, \quad 0 \leqslant t \leqslant T,
\end{equation}
is a solution, at least a viscosity solution, to the PDE (\ref{EDP}).

\paragraph{Assumption (A.3).}
Let $b : [0,T] \times \mathbb{R}^d \rightarrow \mathbb{R}^d$ and $\sigma : [0,T] \rightarrow \mathbb{R}^{d \times d}$ be continuous  functions and let us assume that there exists  $K \geqslant 0$ such that:
\begin{enumerate}
 \item for all $t \in [0,T]$, $\abs{b(t,0)} \leqslant K$, and 
$$\forall (x,x') \in \mathbb{R}^d \times \mathbb{R}^d, \quad \abs{b(t,x)-b(t,x')} \leqslant K \abs{x-x'};$$
 \item $\sigma$ is bounded.
\end{enumerate}
\begin{lem}
\label{moments expo X^2}
 $$\forall \lambda \in \left[0,\frac{1}{2e^{2K T}\norm{\sigma}_{\infty}^2 T} \right[, \exists C_T \geqslant 0, \exists C \geqslant 0, \quad \mathbb{E}\left[ \sup_{0 \leqslant t \leqslant T} e^{\lambda \abs{X_t^{t_0,x_0}}^2} \right] \leqslant C_T e^{C\abs{x_0}^2},$$
with $T \mapsto C_T$ nondecreasing.
\end{lem}
\paragraph{Proof.}
As in \cite{Briand-Hu-08} we easily show that, for all $\eps >0$, we have
\begin{eqnarray*}
\sup_{0 \leqslant t \leqslant T} \abs{X_t^{t_0,x_0}} &\leqslant & \left( \abs{x_0}+K T + \sup_{0 \leqslant t \leqslant T} \abs{\int_0^t \mathbbm{1}_{s \geqslant t_0}\sigma(s) dW_s} \right) e^{K T}\\
\sup_{0 \leqslant t \leqslant T} \abs{X_t^{t_0,x_0}}^2 &\leqslant & C_{\eps}(T^2+ \abs{x_0}^2)+ (1+\eps)e^{2K T} \sup_{0 \leqslant t \leqslant T} \abs{\int_0^t \mathbbm{1}_{s \geqslant t_0}\sigma(s) dW_s}^2. 
\end{eqnarray*}
We define $\tilde{\lambda}:=\lambda(1+\eps)e^{2K T}$. It follows from the Dambis-Dubins-Schwarz representation theorem and the Doob's maximal inequality that 
$$\mathbb{E} \left[ \sup_{0 \leqslant t \leqslant T} \exp \left( \tilde{\lambda} \abs{\int_0^t \mathbbm{1}_{s \geqslant t_0} \sigma(s)dW_s}^2 \right) \right] \leqslant \mathbb{E} \left[ \sup_{0 \leqslant t \leqslant \norm{\sigma}_{\infty}^2T} e^{\tilde{\lambda} \abs{W_t}^2} \right] \leqslant 4\mathbb{E} \left[  e^{\tilde{\lambda}  \norm{\sigma}_{\infty}^2T \abs{W_1}^2} \right],$$
which is a finite constant if $\tilde{\lambda}  \norm{\sigma}_{\infty}^2T < 1/2$. \cqfd
With this observation in hands, we can give our assumptions on the nonlinear term of the PDE and the terminal condition.

\paragraph{Assumption (A.4).}
Let $g: [0,T] \times \mathbb{R}^d \times \mathbb{R} \times \mathbb{R}^d \rightarrow \mathbb{R}$ and $h:\mathbb{R}^d \rightarrow  \mathbb{R}$ be continuous and let us assume moreover that there exist five constants  $r\geqslant 0$, $\beta \geqslant 0$, $\gamma \geqslant 0$, $\alpha \geqslant 0$ and $\alpha'\geqslant 0$ such that:
\begin{enumerate}
\item for each $(t,x,z) \in [0,T] \times \mathbb{R}^d \times \mathbb{R}^{1\times d}$,
$$\forall (y,y') \in \mathbb{R}^2, \quad \abs{g(t,x,y,z)-g(t,x,y',z)} \leqslant \beta \abs{y-y'};$$
\item for each $(t,x,y) \in [0,T] \times \mathbb{R}^d \times \R, \quad z\mapsto g(t,x,y,z)$ is convex on $\R^{1\times d}$;
\item for each $(t,x,y,z) \in [0,T] \times \R^d \times \R \times \R^{1\times d}$,
$$-r(1+ \abs{x}^2 +  \abs{y} + \abs{z}) \leqslant g(t,x,y,z) \leqslant r+ \alpha \abs{x}^2 + \beta \abs{y} +\frac{\gamma}{2} \abs{z}^2,$$
$$-r-\alpha' \abs{x}^2  \leqslant h(x) \leqslant r(1+  \abs{x}^2);$$
\item for each $(t,x,x',y,z) \in [0,T] \times \R^d \times \R^d \times \R \times \R^{1\times d}$,
$$\abs{ g(t,x,y,z)-g(t,x',y,z)} \leqslant r(1+\abs{x}+\abs{x'})\abs{x-x'},$$
$$\abs{h(x)-h(x')} \leqslant r(1+ \abs{x}+\abs{x'})\abs{x-x'};$$
\item $$\alpha'+T\alpha < \frac{1}{2\gamma e^{3\beta T} \norm{\sigma}_{\infty}^2 T}.$$
\end{enumerate}
Thanks to Lemma~\ref{moments expo X^2}, we see that there exist $q >\gamma e^{\beta T}$ and $\eps>0$ such that $h^-(X_T^{t_0,x_0}) + \int_0^T \Big(C + \alpha \abs{X_t^{t_0,x_0}}^2\Big) dt$ has an exponential moment of order $q$ and $h^+(X_T^{t_0,x_0})+\int_0^T \left(r+r \abs{X_t^{t_0,x_0}}^2\right) dt$ has an exponential moment of order $\eps$. So we are able to apply Corollary~\ref{cor existence} and Theorem~\ref{thm unicite} to construct a unique solution $(Y^{t_0,x_0},Z^{t_0,x_0})$ to the BSDE~(\ref{EDSR markovienne}). Let us prove that $u$ is a viscosity solution to the PDE~(\ref{EDP}).
\begin{prop}
\label{sol de viscosite}
 Let assumptions (A.3) and (A.4) hold. The function $u$ defined by~(\ref{feynman-kac}) is continuous on $[0,T] \times \R^d$ and satisfies
$$\forall (t,x) \in [0,T] \times \R^d, \quad \abs{u(t,x)} \leqslant C(1+\abs{x}^2).$$
Moreover $u$ is a viscosity solution to the PDE (\ref{EDP}).
\end{prop}
Before giving a proof of this result, we will show some auxiliary results about admissible control sets. We have already notice in Remark~\ref{N=1} that we have a simpler representation theorem when $T$ is small enough to take $N=1$ in (\ref{definition N}). So we define a constant $T_1>0$ such that for all $T \in [0,T_1]$ we are allowed to set $N=1$. We will reuse notations of section~\ref{uniqueness result section}. By using Remark~\ref{remarque espe q2 bornee}, for all $ T \in [0,T_1]$, $t \in [0,T]$, $x \in \mathbb{R}^d$, we define the admissible control set
\begin{eqnarray*}
 \mathcal{A}_{0,T}(t,x) & := & \left\{(q_s)_{s \in [0,T]}, \quad \int_{0}^{T} \abs{q_s}^2ds < +\infty \,\,\,\, \mathbb{P}-a.s., \quad  \mathbb{E}^{\mathbb{Q}} \left[\int_{0}^{T} \abs{q_s}^2ds \right] < + \infty,\right. \\
& &(M_t)_{t \in [0,T]} \textrm{ is a martingale}, \quad \mathbb{E}^{\mathbb{Q}} \left[\abs{h(X_T^{t,x})}+\int_{0}^{T} \abs{f(s,X_s^{t,x},0,q_s)}ds \right]<+\infty, \\
& & \left. \textrm{ with } M_t:=\exp \left(\int_0^t q_s dW_s-\frac{1}{2}\int_0^t \abs{q_s}^2ds \right) \textrm{ and } \frac{d\mathbb{Q}}{d\mathbb{P}}:=M_T \right\}.
\end{eqnarray*}
We will prove a first lemma and then we will use it to show that this admissible control set does not depend on $t$ and $x$.
\begin{lem}
\label{lemme esp sous Q de X}
$\exists C>0$ such that $\forall  T \in[0,T_1]$, $\forall t \in [0,T]$, $\forall x \in \mathbb{R}^d$, $\forall q \in \mathcal{A}_{0,T}(t,x)$, $\forall s \in [t,T]$,
$$\mathbb{E}^{\mathbb{Q}} \left[\abs{X_s^{t,x}}^2 \right] \leqslant C \left(1+ \abs{x}^2+T \int_t^s \mathbb{E}^{\mathbb{Q}} \left[\abs{q_u}^2 \right]du \right).$$
\end{lem}
\begin{rem}
 $q$ and $\mathbb{Q}$ depend on $x$ and $t$ but we do not write it to simplify notations.
\end{rem}
\paragraph{Proof.}
For all $s \in [t,T]$ we have an obvious inequality
\begin{eqnarray*}
 \abs{X_s^{t,x}}^2 & \leqslant & C \left(  1+\abs{x}^2+\left(\int_t^s \abs{X_u^{t,x}}du\right)^2+\sup_{t \leqslant t' \leqslant T} \abs{\int_t^{t'} \sigma(u)dW_u^q}^2+\left(\int_t^s \abs{q_u}du \right)^2 \right).
\end{eqnarray*}
Then, by applying Cauchy-Schwarz's inequality and Doob's maximal inequality, we obtain
\begin{eqnarray*}
 \mathbb{E}^{\mathbb{Q}} \left[\abs{X_s^{t,x}}^2 \right] & \leqslant & C \left(  1+\abs{x}^2+T \int_t^s \mathbb{E}^{\mathbb{Q}} \left[\abs{X_u^{t,x}}^2 \right]du+\mathbb{E}^{\mathbb{Q}}\left[\abs{\int_t^{T} \sigma(u)dW_u^q}^2\right]\right.\\
& &\left.+T\mathbb{E}^{\mathbb{Q}} \left[\int_t^s \abs{q_u}^2du \right] \right).
\end{eqnarray*}
Finally, the Gronwall's Lemma gives us the result.
\cqfd

\begin{prop}
  $\mathcal{A}_{0,T}(t,x)$ is independent of $t$ and $x$, so we will write it $\mathcal{A}_{0,T}$.
\end{prop}
\paragraph{Proof. }
Let $x,x' \in \R^d$, $t,t' \in [0,T]$ and $q \in \mathcal{A}_{0,T}(t,x)$. We will show that $q \in \mathcal{A}_{0,T}(t',x')$. Firstly,
$$\mathbb{E}^{\mathbb{Q}} \left[\abs{h(X_T^{t',x'})} \right] \leqslant C\left( 1+ \mathbb{E}^{\mathbb{Q}} \left[\abs{X_T^{t',x'}}^2 \right] \right) \leqslant C\left( 1+ \int_{t'}^{T} \mathbb{E}^{\mathbb{Q}} \left[ \abs{q_u}^2 \right] du \right) < +\infty.$$
Moreover
$$-C(1+\abs{X_s^{t',x'}}^2) \leqslant \frac{1}{2\gamma} \abs{q_s}^2-C(1+\abs{X_s^{t',x'}}^2) \leqslant f(s,X_s^{t',x'},0,q_s),$$
and
$$ f(s,X_s^{t',x'},0,q_s) \leqslant f(s,X_s^{t,x},0,q_s)+C(\abs{X_s^{t,x}}^2+\abs{X_s^{t',x'}}^2).$$
So, $\abs{f(s,X_s^{t',x'},0,q_s)} \leqslant \abs{f(s,X_s^{t,x},0,q_s)}+C(\abs{X_s^{t,x}}^2+\abs{X_s^{t',x'}}^2)$ and finally $$\mathbb{E}^{\mathbb{Q}} \left[\int_0^T \abs{f(s,X_s^{t',x'},0,q_s)}ds \right] < +\infty.$$
\cqfd
Now we will do a new restriction of the admissible control set.
\begin{prop}
\label{majoration Y et E[q]}
 $\exists T_2 \in ]0,T_1]$, $\exists \tilde{C}>0$, such that, $\forall T \in[0,T_2]$, $\forall t \in [0,T]$, $\forall s \in [0,T]$, $\forall x \in \mathbb{R}^d$,
$$\abs{Y_s^{t,x}} \leqslant \tilde{C}(1+\abs{x}^2) \quad \textrm{ and } \quad \mathbb{E}^{\mathbb{Q}^*} \left[\int_t^T \abs{q_u^*}^2 du \right] \leqslant \tilde{C}(1+\abs{x}^2).$$
\end{prop}
\paragraph{Proof. }
We are able to use estimations of the existence Theorem~\ref{existence} and Lemma~\ref{moments expo X^2}:
\begin{eqnarray*}
-C\log \mathbb{E} \left[\sup_{0 \leqslant s \leqslant T} \exp \left( C +\gamma e^{\beta T}(\alpha'+T\alpha)\abs{X_s^{t,x}}^2\right) \right] &\leqslant Y_s^{t,x} \leqslant& C\left(1 +\mathbb{E}\left[ \sup_{0 \leqslant s \leqslant T} \abs{X_s^{t,x}}^4 \right] \right)^{1/2}\\
 -\tilde{C}(1+\abs{x}^2) &\leqslant Y_s^{t,x} \leqslant& \tilde{C}(1+\abs{x}^2).\\
\end{eqnarray*}
Then, according to the representation theorem, we have
\begin{eqnarray*}
Y_0^{t,x}& =& \mathbb{E}^{\mathbb{Q}^*} \left[ h(X_T^{t,x})+\int_0^T f(s,X_s^{t,x},Y_s^{t,x},q_s^*)ds \right]\\
 & \geqslant & -C-\alpha' \mathbb{E}^{\mathbb{Q}^*} \left[ \abs{X_T^{t,x}}^2 \right]\\
& & +\frac{1}{2\gamma} \mathbb{E}^{\mathbb{Q}^*} \left[ \int_0^T \abs{q_u^*}^2 du \right] - \alpha \mathbb{E}^{\mathbb{Q}^*} \left[\int_0^T \abs{X_s^{t,x}}^2ds \right] - \beta \mathbb{E}^{\mathbb{Q}^*} \left[ \int_0^T \abs{Y_s^{t,x}}ds \right].\\
\end{eqnarray*}
But, thanks to the uniqueness, we have $Y_s^{t,x}=Y_s^{s,X_s^{t,x}}$ for $s\geqslant t$, so
$\mathbb{E}^{\mathbb{Q}^*} \left[ \abs{Y_s^{t,x}} \right] \leqslant C\left(1+\mathbb{E}^{\mathbb{Q}^*} \left[ \abs{X_s^{t,x}}^2 \right]\right)$. Moreover, we are allowed to use Lemma~\ref{lemme esp sous Q de X},
\begin{eqnarray*}
 Y_0^{t,x}& \geqslant & -C(1+\abs{x}^2) -C(\alpha'+T\alpha+\beta C) \left(1+ \abs{x}^2 +T \int_t^T \mathbb{E}^{\mathbb{Q}^*} \left[ \abs{q_u^*}^2 \right] du \right)\\
& &+\frac{1}{2\gamma} \mathbb{E}^{\mathbb{Q}^*} \left[ \int_0^T \abs{q_u^*}^2 du \right],\\
& \geqslant & -C(1+\abs{x}^2)+\left( \frac{1}{2\gamma} -CT \right) \mathbb{E}^{\mathbb{Q}^*} \left[ \int_0^T \abs{q_u^*}^2 du \right].
\end{eqnarray*}
We set $0<T_2 \leqslant T_1$ such that $\frac{1}{2\gamma} -CT>0$ for all $T \in[0,T_2]$.
Finally,
$$\mathbb{E}^{\mathbb{Q}^*} \left[ \int_0^T \abs{q_u^*}^2 du \right] \leqslant C(1+\abs{x}^2)+Y_0^{t,x} \leqslant \tilde{C}(1+\abs{x}^2).$$
\cqfd
According to the Proposition~\ref{majoration Y et E[q]} we know that $\mathbb{E}^{\mathbb{Q}^*} \left[\int_t^T \abs{q_u^*}^2 du \right] \leqslant \tilde{C}(1+\abs{x}^2)$ so we are allowed to restrict $\mathcal{A}_{0,T}$: for all $r\geqslant 0$ we define 
\begin{equation}
 \mathcal{A}_{0,T}^r=\mathcal{A}_{0,T} \cap \set{(q_s)_{s \in [0,T]}, \mathbb{E}^{\mathbb{Q}} \left[\int_t^T \abs{q_u}^2 du \right] \leqslant \tilde{C}(1+r^2)}.
\end{equation}
With this new admissible control set we will prove a last inequality:
\begin{prop}
\label{majoration delta X}
 $\exists C\geqslant0$, $\forall T \in [0,T_2]$, $\forall t,t' \in [0,T]$, $\forall x,x' \in \mathbb{R}^d$, $\forall q \in \mathcal{A}_{0,T}^{\abs{x} \vee |x'|}$, $\forall s \in [0,T]$,
$$\mathbb{E}^{\mathbb{Q}} \left[ \abs{X_s^{t,x}-X_s^{t',x'}}^2 \right] \leqslant C\left( \abs{x-x'}^2+(1+\abs{x}^2+\abs{x'}^2)\abs{t-t'}\right).$$
\end{prop}
\paragraph{Proof. }
$$\mathbb{E}^{\mathbb{Q}} \left[ \abs{X_s^{t,x}-X_s^{t',x'}}^2 \right] \leqslant \mathbb{E}^{\mathbb{Q}} \left[ \abs{X_s^{t,x}-X_s^{t,x'}}^2 \right] +\mathbb{E}^{\mathbb{Q}} \left[ \abs{X_s^{t,x'}-X_s^{t',x'}}^2 \right].$$
We have, for $s\geqslant t$,
$$X_s^{t,x}-X_s^{t,x'}= x-x'+\int_t^s \left( b(u,X_u^{t,x})-b(u,X_u^{t,x'}) \right)du.$$
So,
$$\mathbb{E}^{\mathbb{Q}} \left[ \abs{X_s^{t,x}-X_s^{t,x'}}^2 \right] \leqslant C\left( \abs{x-x'}^2+\int_t^s \mathbb{E}^{\mathbb{Q}}\left[ \abs{X_u^{t,x}-X_u^{t,x'}}^2 \right]du \right).$$
We apply Gronwall's Lemma to obtain that
$$\mathbb{E}^{\mathbb{Q}} \left[ \abs{X_s^{t,x}-X_s^{t,x'}}^2 \right] \leqslant C\abs{x-x'}^2.$$
Now we deal with the second term. Let us assume that $t \leqslant t'$. For $s \leqslant t$, $X_s^{t,x'}-X_s^{t',x'}=0$. When $t \leqslant s \leqslant t'$, we have
$$X_s^{t,x'}-X_s^{t',x'}=\int_t^s b(u,X_u^{t,x'})du+\int_t^s \sigma(u) dW_u^q+\int_t^s \sigma(u)q_u du.$$
So,
\begin{eqnarray*}
 \mathbb{E}^{\mathbb{Q}} \left[ \abs{X_s^{t,x'}-X_s^{t',x'}}^2 \right] &\leqslant& C\left(   \mathbb{E}^{\mathbb{Q}} \left[ \left(\int_t^{t'}\abs{b(u,X_u^{t,x'})} du\right)^2\right] + \int_t^{t'} \abs{\sigma(u)}^2du+ \mathbb{E}^{\mathbb{Q}} \left[ \left(\int_t^{t'}\abs{\sigma(u)q_u} du\right)^2\right]\right)\\
& \leqslant & C\left( \abs{t'-t}+\abs{t'-t}\int_t^{t'} \mathbb{E}^{\mathbb{Q}} \left[ \abs{X_u^{t,x'}}^2\right]du+ \abs{t'-t} \int_t^{t'} \mathbb{E}^{\mathbb{Q}} \left[ \abs{q_u}^2\right]du \right)\\
& \leqslant & C\abs{t'-t}\left( 1+\abs{x'}^2+\int_0^{T} \mathbb{E}^{\mathbb{Q}} \left[ \abs{q_u}^2\right]du \right)\\
& \leqslant & C(1+\abs{x}^2+\abs{x'}^2)\abs{t'-t}.
\end{eqnarray*}
Lastly, when $t' \leqslant s$,
$$X_s^{t,x'}-X_s^{t',x'}=X_{t'}^{t,x'}-X_{t'}^{t',x'}+\int_{t'}^s b(u,X_u^{t,x'})-b(u, X_u^{t',x'}) du.$$
So,
$$\mathbb{E}^{\mathbb{Q}} \left[ \abs{X_s^{t,x'}-X_s^{t,x'}}^2 \right] \leqslant C(1+\abs{x}^2+\abs{x'}^2)\abs{t'-t} + \int_{t'}^s \mathbb{E}^{\mathbb{Q}} \left[ \abs{X_u^{t,x'}-X_u^{t',x'}}^2 du \right],$$
and according to Gronwall's Lemma,
$$\mathbb{E}^{\mathbb{Q}} \left[ \abs{X_s^{t,x'}-X_s^{t,x'}}^2 \right] \leqslant C(1+\abs{x}^2+\abs{x'}^2)\abs{t'-t}.$$
\cqfd

\paragraph{Proof of Proposition~\ref{sol de viscosite}.}
First of all, let us assume that $T<T_2$. With this condition, we are allowed to use all previous propositions. Firstly, the quadratic increase of $u$ is already proved in Proposition~\ref{majoration Y et E[q]}. Then, we will show continuity of $u$ in $(t_0,x_0) \in [0,T] \times \mathbb{R}^d$. We have
$$\forall (t,x) \in [0,T] \times \mathbb{R}^d, \quad \abs{u(t,x)-u(t_0,x_0)} \leqslant \abs{u(t,x)-u(t,x_0)}+\abs{u(t,x_0)-u(t_0,x_0)}.$$
Let us begin with the fist term. We define $r:=\abs{x}\vee \abs{x_0}$.Thanks to the representation theorem, we have
$$Y_t^{t,x} = \essinf_{q \in \mathcal{A}_{0,T}^r} Y_t^{q,t,x} \quad \textrm{and} \quad  Y_t^{t,x_0} = \essinf_{q \in \mathcal{A}_{0,T}^r} Y_t^{q,t,x_0}.$$
So, 
$$\abs{Y_t^{t,x} - Y_t^{t,x_0}} \leqslant \esssup_{q \in \mathcal{A}_{0,T}^r} \abs{Y_t^{q,t,x} - Y_t^{q,t,x_0}}.$$
But, for $t \leqslant s \leqslant T$,
\begin{eqnarray*}
\abs{Y_s^{q,t,x}-Y_s^{q,t,x_0}} &=& \Big|\mathbb{E}^{\Q}\Big[ h(X_T^{t,x})-h(X_T^{t,x_0}) \\
& &  +\int_s^T \left(f(u,X_u^{t,x},Y_u^{q,t,x},q_u)-f(u,X_u^{t,x_0},Y_u^{q,t,x_0},q_u)\right)du \Big| \mathcal{F}_s\Big]\Big|\\
 & \leqslant & \mathbb{E}^{\Q}\left[C(1+\abs{X_T^{t,x}}^2+\abs{X_T^{t,x_0}}^2) \right]^{1/2} \mathbb{E}^{\Q}\left[\abs{X_T^{t,x}-X_T^{t,x_0}}^2 \right]^{1/2}\\
& & +\int_s^T \mathbb{E}^{\Q}\left[C(1+\abs{X_u^{t,x}}^2+\abs{X_u^{t,x_0}}^2)\right]^{1/2} \mathbb{E}^{\Q}\left[\abs{X_u^{t,x}-X_u^{t,x_0}}^2 \right]^{1/2} du\\
& & + C\int_s^T \mathbb{E}^{\Q}\left[\abs{Y_u^{q,t,x}-Y_u^{q,t,x_0}} \right]du,\\
\end{eqnarray*}
thanks to Assumption (A.4) and Hölder's inequality.
According to Lemma~\ref{lemme esp sous Q de X}, the definition of $\mathcal{A}_{0,T}^r$ and Proposition~\ref{majoration delta X}, we obtain
$$\mathbb{E}^{\Q} \left[\abs{Y_s^{q,t,x}-Y_s^{q,t,x_0}}\right] \leqslant C(1+\abs{x}^2+\abs{x_0}^2)^{1/2}\abs{x-x_0} + C\int_s^T \mathbb{E}^{\Q}\left[\abs{Y_u^{q,t,x}-Y_u^{q,t,x_0}} \right]du.$$
Then, Gronwall's lemma gives us 
$\abs{Y_t^{q,t,x}-Y_t^{q,t,x_0}} \leqslant C(1+\abs{x}+\abs{x_0})\abs{x-x_0}$. Since this bound is independent of $q$, we finally obtain that
$$ \abs{Y_t^{t,x} - Y_t^{t,x_0}} \leqslant C(1+\abs{x}+\abs{x_0})\abs{x-x_0}.$$
Now, we will study the second term. Without loss of generality, let us assume that $t <t_0$.
\begin{eqnarray*}
\abs{Y_t^{t,x_0}-Y_{t_0}^{t_0,x_0}} & \leqslant & \abs{Y_t^{t,x_0}-Y_t^{t_0,x_0}}+\int_t^{t_0} \abs{g(s,x_0,Y_s^{t_0,x_0},0)} ds,\\
& \leqslant & \abs{Y_t^{t,x_0}-Y_t^{t_0,x_0}} + \int_t^{t_0} C(1+\abs{x_0}^2 + \abs{Y_s^{t_0,x_0}}) ds.
\end{eqnarray*}
We apply Proposition~\ref{majoration Y et E[q]} to obtain
$$\abs{Y_t^{t,x_0}-Y_{t_0}^{t_0,x_0}}  \leqslant  \abs{Y_t^{t,x_0}-Y_t^{t_0,x_0}} +C(1+\abs{x_0}^2)(t-t_0).$$
We still have
$$\abs{Y_t^{t,x_0} - Y_t^{t_0,x_0}} \leqslant \esssup_{q \in \mathcal{A}_{0,T}^r} \abs{Y_t^{q,t,x_0} - Y_t^{q,t_0,x_0}}.$$
Moreover, exactly as the bound estimation for $\mathbb{E}^{\Q} \abs{Y_s^{q,t,x}-Y_s^{q,t,x_0}}$, we have, for $t \leqslant s \leqslant T$,
\begin{eqnarray*}
\mathbb{E}^{\Q}\left[ \abs{Y_s^{q,t,x_0}-Y_s^{q,t_0,x_0}} \right] & \leqslant & \mathbb{E}^{\Q}\left[C(1+\abs{X_T^{t,x_0}}^2+\abs{X_T^{t_0,x_0}}^2) \right]^{1/2} \mathbb{E}^{\Q}\left[\abs{X_T^{t,x_0}-X_T^{t_0,x_0}}^2 \right]^{1/2}\\
& & +\int_s^T \mathbb{E}^{\Q}\left[C(1+\abs{X_u^{t,x_0}}^2+\abs{X_u^{t_0,x_0}}^2)\right]^{1/2} \mathbb{E}^{\Q}\left[\abs{X_u^{t,x_0}-X_u^{t_0,x_0}}^2 \right]^{1/2} du\\
& & + C\int_s^T \mathbb{E}^{\Q}\left[\abs{Y_u^{q,t,x_0}-Y_u^{q,t_0,x_0}} \right]du.\\
\end{eqnarray*}
According to Lemma~\ref{lemme esp sous Q de X}, the definition of $\mathcal{A}_{0,T}^r$, Proposition~\ref{majoration delta X} and Gronwall's Lemma, we obtain
$\abs{Y_t^{q,t,x_0}-Y_t^{q,t_0,x_0}} \leqslant C(1+\abs{x}^2+\abs{x_0}^2)\abs{t-t_0}^{1/2}$. Since this bound is independent of $q$, we finally obtain that
$$ \abs{Y_t^{t,x} - Y_t^{t,x_0}} \leqslant C(1+\abs{x}^2+\abs{x_0}^2)\abs{t-t_0}^{1/2}.$$
So,
$$ \abs{u(t,x)-u(t_0,x_0)} \leqslant C(1+\abs{x}+\abs{x_0})\abs{x-x_0} + C(1+\abs{x}^2+\abs{x_0}^2)\abs{t-t_0}^{1/2}.$$
We now return to the general case (for $T$) : we set $N \in \mathbb{N}$ such that $T/N < T_2$ and, for $i \in \set{0,...,N}$, we define $t_i:=iT/N$. According to the beginning of the proof, $u$ is continuous on $[t_{N-1},T] \times \mathbb{R}^d$. We define $h_{N-1}(x):=Y_{t_{N-1}}^{t_{N-1},x}$. Since $\abs{h_{N-1}(x)-h_{N-1}(x')} \leqslant C(1+\abs{x}+\abs{x'})\abs{x-x'}$, we are allowed to reuse previous results to show the continuity of $u$ on $[t_{N-2},T_{N-1}] \times \mathbb{R}^d$. Thus, we can iterate this argument to show the continuity of $u$ on $[0,T] \times \R^d$. Moreover the quadratic increase of $u$ with respect to the variable $x$ results from the quadratic increase of $u$ on each interval.

Finally, we will use a stability result to show that $u$ is a viscosity solution to the PDE (\ref{EDP}). As in the proof of Theorem~\ref{existence}, let us consider the function
$$g_n(t,x,y,z):=\inf \set{g(t,x,p,q)+n\abs{p-y}+n\abs{q-z}, (p,q) \in \mathbb{Q}^{1+d}}.$$
We have already seen that $(g_n)_{n \geqslant \lceil r\rceil}$ is increasing and converges uniformly on compact sets to $g$. Let $(Y^{n,t,x},Z^{n,t,x})$ be the unique solution in $\mathcal{S}^2 \times M^2$ to BSDE($h(X_T^{t,x})$,$-g_n(.,X_.^{t,x},.,.)$). We define $u_n(t,x):=Y^{n,t,x}_t$. Then by a classical theorem (see e.g. \cite{ElKaroui-Peng-Quenez-97}), $u_n$ is a viscosity solution to the PDE
\begin{equation*}
\partial_t u(t,x)+\mathcal{L}u(t,x)-g_n(t,x,u(t,x),-\sigma^*\nabla_xu(t,x))=0, \quad u(T,.)=h.
\end{equation*}
Moreover, it follows from the classical comparison theorem that $(u_n)_{n \geqslant \lceil r\rceil}$ is decreasing and, by construction, converges pointwise to $u$. Since $u$ is continuous, Dini's theorem implies that the convergence is also uniform on compacts sets. Then, we apply a stability result (see e.g. Theorem~1.7. of chapter~5 in \cite{Bardi-CapuzzoDolcetta-97}) to prove that $u$ is a viscosity solution to the PDE (\ref{EDP}).
\cqfd

{\bf Remark.} The uniqueness of viscosity solution to PDE is considered by Da Lio and Ley in \cite{DaLio-Ley-06} and \cite{DaLio-Ley-08}.

\end{document}